\documentclass[sustainability,article,accept,pdflatex,moreauthors,a4paper]{mdpi} 
\usepackage{tikz}
\usetikzlibrary{arrows.meta}
\usetikzlibrary{fadings}
\usepgflibrary{shadings}
\tikzfading[name=fade inside,
            inner color=transparent!30,
            outer color=transparent!30]

\firstpage{1} 
\pubvolume{14}
\issuenum{1}
\articlenumber{12500}
\pubyear{2022}
\copyrightyear{2022}
\externaleditor{Academic Editor: Jin Su Jeong and David González-Gómez} 
\datereceived{30 July 2022}
\daterevised{} 
\dateaccepted{17 September 2022} 
\datepublished{30 September 2022} 
\hreflink{https://doi.org/10.3390/su141912500} 

\Title{Sustainable Learning, Cognitive Gains, and Improved Attitudes in College Algebra Flipped Classrooms}

\TitleCitation{Sustainable Learning, Cognitive Gains, and Improved Attitudes in College Algebra Flipped Classrooms}

\Author{Natanael Karjanto $^{1,}$*\orcidA{} and Maxima Joyosa Acelajado $^{2}$}

\AuthorNames{N. Karjanto and M.J. Acelajado}

\AuthorCitation{Karjanto, N.; Acelajado, M.J.}

\address{%
$^{1}$ \quad {Department} of Mathematics, University College, Natural Science Campus Sungkyunkwan University, Suwon~16419, Korea\\
$^{2}$ \quad {Department} of Science Education, College of Education De La Salle University, 2401 Taft Ave, Malate, Manila~1004,  Philippines}

\corres{Correspondence: natanael@skku.edu}

\abstract{To respond to global issues positively, education systems in higher education institutions play a significant role in empowering learners as well as promoting sustainable development goals. By implementing curricula that cultivate cross-cutting and transversal key competencies for sustainability, such as critical thinking, problem-solving, and collaboration, we prepare our pupils to become sustainability citizens, who not only sustain learning throughout their lives in various circumstances and across different disciplines but also engage constructively and responsibly toward any future world's challenges through their dispositions, strategies, and skills. One such sustainable teaching methodology is known as the flipped classroom, an active-learning, student-centered, flexible, and multidimensional pedagogy. Our objective is to investigate the effect of such pedagogy on learners' academic achievement and their attitude toward mathematics using both quantitative and qualitative methods. We cultivated sustainable learning in mathematics education for college freshmen ($n = 55$) by exposing them to both the conventional teaching method (CTM) and flipped classroom pedagogy (FCP). By splitting them into control and experimental groups alternately ($n_1 = 24$, $n_2 = 31$) and by selecting the four most challenging topics in college algebra, we measured their cognitive gains quantitatively via a sequence of pre- and post-tests. The topics are factorization, rational expressions, radical operations, and applied problems. Both groups improved academically over time across all these four topics with statistically very significant outcomes $(p < 0.001)$. Although they were not always statistically significant ($p > 0.05$) in some topics, the post-test results suggest that generally, the FCP trumps the CTM in cognitive gains, except for the first topic on factorization, where the opposite is true with a very statistically significant mean difference $(p < 0.001)$. By examining non-cognitive gains qualitatively, we analyzed the students' feedback on the FCP and their responses to a perception inventory. The finding suggests a favorable response toward the FCP with primary improvements in the attitudes toward mathematics and increased levels of cooperation among students. Since these students are so happy to have control of their own learning, they were more relaxed, motivated, confident, active, and responsible in learning under the FCP. We are confident that although this study is relatively small in scale, it will yield incremental and long-lasting effects not only for the learners themselves but also for other role-takers in education sectors who aspire in nurturing sustainable long-life learning and achieving sustainable development goals successfully.}


\begin{document}

\section{Introduction} \vspace*{1mm}

One of the 17 United Nations' Sustainable Development Goals (SDGs) for achieving a better and more sustainable future for all beings is ``quality education'' (Goal 4, herein SDG4). The two facets of SDG4 ensure inclusive and equitable quality education and promote lifelong learning opportunities for all~\cite{un2030sdg}. Although the primary targets of SDG4 are children and youths in their primary and secondary educations, the principle can be extended into higher education institutions (HEIs) as well. HEIs play an essential role in promoting sustainability and advocating the SDGs by educating and training prospective teachers, policymakers, future leaders, entrepreneurs, and other professionals~\cite{hesi2021higher}.  

In particular, through their cutting-edge research and continuously improving their curricula, HEIs advance sustainability competencies among their pupils. The specialized agency UNESCO also supports SDG4 through education for sustainable development (ESD) and global citizenship education (GCED). In turn, this education will enhance necessary cross-cutting key competencies for sustainability that are relevant not only to SDG4 but also to all SDGs. Among these crucial key competencies are critical thinking, integrated problem-solving, and collaboration competencies~\cite{unesco2007education}. 

Keeping this association in mind, having a solid foundation in science, technology, engineering, art, and mathematics (STEAM) education is very essential in order to develop sustainable citizens since these subjects inherently cultivate critical thinking and problem-solving skills, the very qualities which are aligned with ESD and GCED in SDG4~\cite{rogovaya2019critical}. These skills can be used throughout life irrespective of whatever professions the students might choose to pursue. The very basics of all these disciplines, namely mathematics as well as the way we teach and how students learn about it (mathematics education), should be of particular interest to many educators around the world who aspire to instill lifelong learning among their students, and at the same time, nourish them with the ESD and GCED transversal competencies.

Thanks to the progress in information technology, many novel pedagogical methodologies have emerged and received more attention during the past two decades, particularly those that emphasized student-centered learning activities. These include but are not limited to flipped classroom, blended learning, problem and/or project-based learning, inquiry-based learning, collaborative learning, and inclusive education, including its framework of universal learning design, among others. The body of literature confirms that these educational pedagogies demonstrated excellent educational outcomes and provided platforms for supportive learning environments that align with the SDGs~\cite{kioupi2019education,galan2020sustainable}. Focusing on a particular pedagogy and observing how it affects learners' cognitive gains and improves attitudes toward a particular subject or learning in more general would be an essential step in cultivating transferable competencies for sustainability. 

In this article, we discuss how flipped classrooms in a college algebra module could improve students' cognitive gains as well as their attitudes toward mathematics. Unless specified differently, throughout this paper, cognitive gains refer to the positive gains in knowledge and comprehension, which can be translated through academic performance. Although non-cognitive gains encompass a wide range of abilities such as communication, teamwork, perseverance, conscientiousness, and motivation, among others, what we refer to is related to the latter, i.e., students' attitudes toward mathematics. Furthermore, we will refer to the two pedagogical approaches used in this study as the FCP and the CTM, which refer to the flipped classroom pedagogy and conventional teaching method, respectively.

With the nature of present-day students and their exposure to a variety of technological tools, nowadays teaching requires a degree of flexibility in addressing diverse learning styles, a wide range of level capability, and handling large classes. A breed of confident and competent problem-solvers who are eager to learn new things on their own and the improvement of students' achievement in and attitudes toward mathematics are very much desired not only to survive in this modern society but also to develop sustainable citizens. Hence, there is a pressing need for advancing from teacher-centered, passive-receptive learning of the CTM into a student-centered, active-learning FCP. Some evidence suggests that the FCP could be a sustainable active-learning pedagogy when learning is disrupted such as during the recent COVID-19 pandemic~\cite{collado2021flipped}. 

This study investigated the impact of the FCP on students' academic achievements in a mathematics course among freshmen majoring in mathematics education at a private, Catholic coeducational research university in the Philippines. Additionally, we are addressing the following research questions:
\begin{itemize}
\item Between the experimental and control groups and among the four most challenging topics in college algebra, which group obtains better cognitive gains? Which topics would be learned better using the FCP?
\item For non-cognitive gains, what are students' reactions to the use of the FCP? What are their attitudes toward mathematics after learning the course via the FCP?
\end{itemize}

The article is organized as follows. After this introduction, Section~\ref{litrev} provides a literature review on sustainable learning in education (SLE) and how the FCP has been successfully implemented in a wide range of disciplines from mathematics to ESD. It then offers a conceptual framework for this study that SDG4 covers both ESD and SLE, and the FCP is one possible option for SLE. Section~\ref{method} covers research objectives and methodology, including the participants and applied measurement. Section~\ref{result} outlines the result of our experiments and discusses what these findings mean. Finally, Section~\ref{conclude} concludes our discussion and provides further recommendations.

\section{Literature Review}		\label{litrev} 	\vspace*{2mm}

\subsection{Sustainable Learning in Mathematics Education}		\label{slime}	\vspace*{1mm}

Ben-Eliyahu (2021) clarified a potentially baffling concept between ``learning sustainable development'' and ``sustainable learning in education (SLE)''~\cite{ben2021sustainable}. The designated terminology for the former as ``sustainable learning'' might be easily mistaken for the latter or vice versa. On the one hand, a better terminology for learning sustainable development would be ``sustainability learning'' or ``ESD'', which refers to an approach to education that emphasize the importance of humans living in harmony with nature, either by specifically teaching sustainability principles or by integrating them into a curriculum by including key issues of sustainable development, as we have mentioned briefly in the introduction~\cite{venkataraman2009education,hopkins2002education,boeve2015effectiveness,huckle2015UN}.

On the other hand, SLE refers to ``learning that lasts''; it can be achieved through well-structured and responsive teaching that matters for all learners~\cite{graham2015sustainable}. It consists of two facets of learning that can be likened to two sides of a coin. On the one side of the coin, the learning that takes place in formal education has lasting value to learners into the future. On the other side of the coin and at the same time, this sustainable type of learning will also encourage pupils to continue their education journey by embracing lifelong learning themselves. Indeed, SLE interlaces marvelously with SDG4; see~\cite{un2030sdg}. Furthermore, as suggested by Hays and Reinders (2020), any type of sustainable learning and education should instill in students the skills and dispositions to thrive in a complicated and challenging world as well as the desire to contribute positively in creating the world a better place to live~\cite{hays2020sustainable}.

Among one of the seven recommendations for creating sustainability education at HEIs, Moore (2005) proposed establishing a space for pedagogical transformation. The endorsement encompasses not only improving the interaction between students and instructors but also promoting student-centered, reflective, critical, transformative, and experiential learning~\cite{moore2005seven}. This advocate is underpinned by a recent article on the learning environment in the context of SDG4, where it is confirmed that an educational environment that accommodates sustainable learning is one that encouraged active role from both the pupils and educators, such as class participation, critical thinking, nurturing curiosity, and cultivating creativity, among others~\cite{galan2020sustainable}. 

Similar to renewable and sustainable energy, SLE preserves the learning process throughout one's life as the scene of the world is changing. It endows learners with skills and strategies to rejuvenate themselves through inquiry, self-assessment, and evaluation of their environment and social systems. SLE encompasses four aspects of self-regulated learning models that exhibit an analogy with sustainable nature: renewing and relearning; independent and collaborative learning; active learning;  transferability~\cite{ben2021sustainable,zimmerman1990self,pintrich1995understanding,boekaerts1999self,puustinen2001models,zimmerman2011self}. As educators, in addition to grooming our students with future-focused experiences and skills, we also ought to nurture their confidence and refine their awareness of achieving positive changes.

On a larger scale, HEIs play a key role in sustainability; again like the two sides of a coin, by promoting both ESD as well as SLE. On the one side of the coin, by facilitating and designing curricula that center around sustainability, i.e., ESD, HEIs contribute crucially to creating a sustainability mindset among their members, including faculty, administrative personnel, and, in particular, the student body, where a new generation of future leaders would emerge~\cite{cortese2003critical,jones2010sustainability,zaleniene2021higher,veidemane2022education}. On the other side of the coin, by improving their didactics and pedagogy, i.e., SLE, HEIs prepare current and prospective learners with transversal competencies that are necessary for tackling not only personal, a smaller-level, challenging circumstances but also global, a larger-level, economic, social, and environmental challenges~\cite{sterling2001sustainable,branden2012sustainable,fullan2015new,peris2015sustainable,sa2018transversal,membrillo2021sustainability}. 

Although the aforementioned cited works administer the field of education in the general context, the same principles are also certainly appropriate for mathematics education. Renert (2011) attempted to address an inquiry about how to reconcile the urgent need to act for a sustainable future with the current practices of mathematics education by presenting a model of possible responses to sustainability in mathematics education~\cite{renert2011mathematics}. The proposed model adapts two existing stage models of approaches to sustainability to the context of mathematics education, i.e., Sterling's (2001) and Edwards' (2010) models of educational responses and organizational approaches to sustainability, respectively~\cite{sterling2001sustainable,edwards2010organizational}. The three types of the educational response of accommodation, reformation, and transformation mean education \emph{about}, \emph{for}, and \emph{as} sustainability, respectively. The former two fit well with ESD, whereas the final type aligns with SLE. 

With an exception of Summer's (2020) exploration of how a sustainable primary mathematics education ought to be implemented~\cite{summer2020sustainable}, other works that intersect between sustainability and mathematics education usually concentrate around ESD instead of SLE, e.g.,~\cite{hamilton2014sustainability,barwell2018some,widiati2019philosophy,li2020philosophy,li2021education,moreno2022training}. By suggesting concrete pedagogical initiatives to tackle primary students' challenges in learning mathematics, Summer (2020) demonstrated that quality education and competent teachers not only decrease their learning difficulties but also equip children with essential mathematical and critical thinking skills~\cite{summer2020sustainable}. These skills, in turn, will eventually provide a solid foundation for acquiring transversal competencies when they grow up, join the workforce, and assess sustainability principles.

\subsection{Flipped Classroom Pedagogy}		\label{fcp}		\vspace*{1mm}

To better comprehend the FCP, we need to agree with what we mean by the traditional teaching and learning approach, dubbed in this article the CTM. In a CTM situation, the students listen to the teacher's lecture-discussion on the day's lesson and they are given homework or assignment to measure their understanding of the lesson. The homework may be discussed in class or simply submitted to the teacher. Teachers are the instruments by which knowledge is communicated. The CTM is primarily teacher-centered where all students are taught the same materials at the same time. The CTM emphasizes direct instruction, predominantly lectures, and a fixed seatwork so that students learn through listening and observation. Classroom instruction is often solely based on textbooks, lectures, and individual written assignments.

An FCP, which is otherwise known as an inverted classroom~\cite{lage2000inverting}, flipped teaching, flipped learning, or the Thayer method~\cite{connors2000thayer,shell2002thayer}, is a combination of viewing video recordings and reading module materials related to the lesson anywhere and anytime before class and applying what has been learned during the face-to-face time~\cite{overmyer2012flipped}. Class time is used more interactively for group discussion, discovery activities, experiments, and class presentations where the teacher's role is to facilitate and assist the students in their quest for further understanding of the lesson. As a result, the students become active learners rather than plain receptacles of information~\cite{king1993from}. 

The FCP is a form of blended learning in which students learn new contents from various modes, such as by viewing video lectures online, browsing websites, reading textbooks, and viewing PowerPoint-like (or Beamer) lecture slide presentations that are either provided by the teacher or the result of their own search mechanisms, usually at home before the class time. What used to be the homework/assignment is now done in class with teachers offering more personalized guidance and interaction with students, instead of lecturing full-time. During the class, the students apply the knowledge they have learned by solving problems, doing practical work, and collaborating with their peers. 

Both classroom and online learning should provide materials that would enable the students to practice and interact with others and still deliver favorable outcomes. Doing the homework to keep them engaged for a deeper understanding of concepts and mastery of skills is done during face-to-face class time. The students work on their own phase and work level in a result-based environment while using online content that serves as the preliminary source of learning and this may be used repeatedly in class or as review materials. The content acquisition comes ahead and concept engagement takes place in class with the students doing interactive activities. By inquiring the teachers on materials that they could not grasp or more challenging topics, positive response and feedback from the teacher could assist both adept and slow learners.

There are several implications in the context of educational theory when implementing the FCP. First, the knowledge becomes personal since each learner might differ in interpretations and thus possesses a distinctive point of view when following a flexible and multi-dimensional pedagogy such as the FCP~\cite{fox2001constructivism}. Second, we allow students to construct knowledge by themselves when they experience different things, rather than passively absorbing it in the case of  the CTM~\cite{arends1998resource,philips1995good}. Third, and as a consequence, learning becomes an active process. Learners construct meaning in their understanding not only through active engagement with their environment but also by establishing meaningful connections between prior and new knowledge~\cite{good1994looking,cooperstein2004beyond}. These aspects constitute a constructivist approach in education theory and we could confirm from the following literature study that the FCP does belong to this approach~\cite{elliott2000educational,vigotsky1978mind}.

According to Strayer (2007), an FCP is a more active, student-centered style of teaching through the use of group projects, discovery activities, experiments, and class presentations that are implemented during classroom time with information-rich, lecture-based direct instruction being used during an out-of-class time, usually delivered through online videos that students view before arriving in class~\cite{strayer2007effects}. In short, this model aims to move the easier parts of teaching and learning into independent practice ahead of learning the more difficult concepts, which are taught face-to-face~\cite{bretzmann2013how}.

Flipped teaching is a pedagogical approach in which direct instruction is moved from the group learning space to the individual learning space, and the resulting group space is transformed into a dynamic, interactive learning environment where the teacher guides students as they apply concepts and engage creatively in the subject matter. The FCP has benefited students who missed classes since they may use the online materials to review and reinforce lessons and it gives them the opportunity to radically rethink how they should use class time effectively~\cite{bergmann2012flip,tucker2012flipped}.

In the FCP, instruction is delivered online and outside the classroom through video, podcasts, or the online learning environment. Homework is moved into the classroom. In this approach, students can take in the information at their own pace and discuss it with the teacher and peers. This creates time in the classroom for collaborative work by the students and more room for a differential approach and remediation by the teacher. Fulton (2012) enumerated the following advantages of the flipped classroom: (1) students move at their own pace; (2) doing ``homework'' in class gives teachers better insight into students' difficulties and learning styles; (3) teachers can more easily customize and update the curriculum and provide it to students anytime; (4) classroom time can be used more effectively and creatively; (5) teachers using the method report seeing increased levels of student achievement, interest, and engagement; (6) learning theory supports the new approaches; (7) the use of technology is flexible and appropriate \mbox{for the 21st-century~learning}~\cite{fulton2012upside}.

Moreover, Herreid and Schiller (2012) argued that with the FCP: (1) there is more time to spend with students on authentic research; (2) students get more time working with scientific equipment that is only available in the classroom; (3) students who missed class can watch the lectures while on the road; (4) the method promotes thinking inside and outside of the classroom; (5) students are more actively involved in the learning process~\cite{herreid2012case}.  The survey conducted by Bishop and Verleger (2013) revealed that most studies on the FCP used single-group study designs to explore student perceptions and reports are somewhat mixed but generally positive overall. Some students tend to prefer in-person lectures to video lectures, but prefer interactive classroom activities over lectures. Anecdotal evidence suggests that student learning is improved for the FCP compared to the CTM~\cite{bishop2013flipped}. 

Hantla (2014) indicated that the flipped classroom is a new iteration of an old way of teaching that enables instructors to do more during face-to-face classroom time than is otherwise possible. The FCP provides an incentive to students to come to class prepared and assess their understanding, focuses on higher-level cognitive activities under the guidance of the teacher, and provides students adequate time to carry out their assignments and get on-the-spot feedback about their work. The FCP is carried out in a learning environment with the support of educational technology where students learn through activity-oriented activities~\cite{hantla2014effects}.

When it comes to teaching mathematics at all levels of education, the body of published literature does not lack examples. The FCP has been successfully implemented in teaching mathematics at the primary~\cite{bergman2016flipped,alrouqi2019flipped,panahi2019reviewing,khalel2021employing} and secondary levels~\cite{saunders2014flipped,bhagat2016impact,makinde2017flipped,lo2017using,lo2017critical,weinhandl2018technology,wei2020effect}. The FCP also encompasses a broad mathematical subjects, such as college algebra~\cite{overmyer2014flipped,vansickle2015adventures}, precalculus~\cite{fulton2012upside}, calculus~\cite{jungic2015flipping,maciejewski2016flipping,mcgivney2013flipping,sahin2015flipping,sonner2015impact,karjanto2019english}, vector calculus~\cite{ziegelmeier2015flipped}, linear algebra~\cite{talbert2014inverting,murphy2016student,novak2017flip,karjanto2017flipped}, statistics~\cite{wilson2013flipped}, and actuarial science~\cite{butt2014student}. Generally, both students' cognitive gains in terms of understanding and academic performance, as well as students' positive attitudes in terms of enjoyment and confidence are affirmative; see also~\cite{poffenberger1959factors,hodges2013improving,lo2017toward,karjanto2017attitude,turra2019flipped}.

Certainly, the successful---and fruitless---attempts of the FCP are not only limited to mathematics education. Many works provide a narrative on pedagogy in other fields, including both STEM and non-STEM fields. In particular, in what follows we provide some evidence in the learning about sustainability and sustainable developments. Buil-Fabreg\'{a} et al. (2019) demonstrated that the FCP has successfully assisted students in improving their transversal competencies and being more conscious of sustainable development requirements~\cite{buil2019flipped}. Rodríguez-Chueca et al. (2020) measured the efficiency of the FCP and challenged-based learning to facilitate learning of sustainability principles and circular economy and discovered that the former is more satisfactory than the latter~\cite{rodriguez2020understanding}. Howell (2021) revealed positive student perceptions when they were exposed to the FCP when enrolling in education for sustainable development courses~\cite{howell2021engaging}.

\subsection{Conceptual Framework}		\vspace*{1mm}

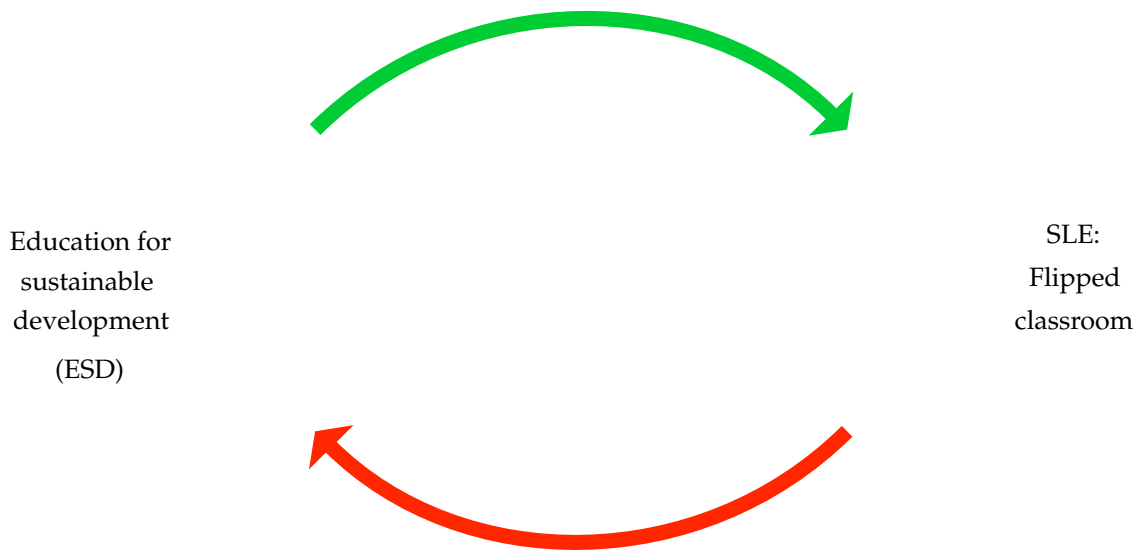
\begin{figure}[H] \vspace*{-0.5cm}
\begin{adjustwidth}{-\extralength}{0cm}
\begin{center}
\begin{tikzpicture}
\tikzset{arrow/.style = {line width=2mm},-{Triangle[length=3mm]}}
\shade[ball color = blue!80!white, path fading = fade inside] (0,0) circle (1.5);
\begin{scope}[xshift=7.3cm]
\shade[ball color = red!80!orange, path fading = fade inside] (0,0) circle (1.5);
\end{scope}
\draw[arrow,draw=green!80!blue] (0,2) to [bend left=45] (7,2);
\draw[arrow,draw=red!70!orange] (7,-2) to [bend left=45] (0,-2);
\node at (-10:-3){Education for};
\node at (0:-3){sustainable};
\node at (11:-3){development};
\node at (22:-3.2){(ESD)};
\node at (3.5:10){SLE:};
\node at (0:10){Flipped};
\node at (-3:10){classroom};
\end{tikzpicture}
\end{center}
\end{adjustwidth}
\caption{A conceptual framework for sustainable development in mathematics education employing SLE and FCP. HEIs play an active role in promoting SDG4 by equipping learners with sustainable learning in (mathematics) education to cope with future challenging circumstances. In particular, the FCP offers quality education and SLE through its flexible teaching yet effective learning approach (top green arrow). Conversely, the FCP can be implemented to cultivate learners in contributing to SDGs via ESD with the aim of embracing and practicing sustainability principles (bottom red arrow).}		\label{theoret}
\end{figure}

Based on the literature review of both SLE as well as the superiority of the FCP over the CTM, we propose a conceptual framework to examine that our study on flipped classrooms in mathematics education can also be expanded to other topics outside mathematics, duplicated, and even improved for better learning outcomes as well as with more positive psychological well-being. Figure~\ref{theoret} displays a conceptual framework for this study.

Our study is based upon a conceptual framework that sustainability and education influence each other and contribute to each other's development. As we reviewed in Subsection~\ref{slime}, SLE in general and sustainable learning in mathematics education in particular require not only strategic HEI curricula that empower learners but also flexible and multi-dimensional teaching approaches that equip learners with transferable skills to renew themselves when facing personal and world challenges~\cite{ben2021sustainable,hays2020sustainable}. Among the highlighted pedagogies that promote SLE are project-based learning and the FCP. Our focus is the latter. Conversely, the body of published literature reviewed in Subsection~\ref{fcp} has demonstrated that the FCP was not only viewed very favorably among students but also compels them to commit to sustainability principles and practice sustainable development when they join the labor market, as demonstrated in several recent studies~\cite{buil2019flipped,rodriguez2020understanding,howell2021engaging}. Hence, as illustrated in Figure~\ref{theoret}, ESD and SLE could function both ways.

\section{Research Material and Methodology}	\label{method}		\vspace*{2mm}

\subsection{Participant}		\vspace*{1mm}

Two classes of 24 and 31 freshmen mathematics education majors at De La Salle University, Manila, the Philippines, enrolled in a three-credit course of College Algebra (course code: TCHALGE) during the first term of the academic year 2014/2015 (September--December~2014), served as respondents of this study. The classes were referred to as Group~1 and Group~2, respectively. The selection methodology was based on the convenience technique since these two classes were simply assigned to and taught by one of us (the second author). The age range of the participants is typically between 18 and 19 years old. The percentages of the female and male students are around 60\% and 40\%, respectively.

\subsection{Measurement}		\vspace*{1mm}

A quasi-experimental design with switching replication was implemented for both groups. The group exposed to the FCP is considered the experimental group, whereas the group exposed to the CTM is regarded as the control group. For the first and second topics, Group~2 was considered the experimental group, while for the third and fourth topics, Group~1 was designated as the experimental group. Both groups were also administered pre- and post-tests on all four topics considered in this study. The two groups were alternately exposed to the FCP and the CTM in the delivery of topics in college algebra which are identified as relatively hard by students who took the subject during preceding terms. These are factoring/factorization, rational expressions, operations on radicals, and solving applied problems. 

On the one hand, the experimental group was provided with a list of all websites, videos, lecture notes, PowerPoint presentation slides, and course modules that they need to watch and read as well as the corresponding assessment materials that they need to respond while viewing the videos or listening to the lectures before class to test their understanding of the topics. During the face-to-face class time, they were given more interactive activities which were done either individually or in groups with the teacher acting as a facilitator. On the other hand, the control group was taught identical topics and was provided with the same interactive activities in class after which the students were given the assignment/homework for submission during the next meeting. 

In our empirical study, we have employed both quantitative and qualitative methods to investigate the effect of the FCP on students' cognitive and non-cognitive gains. The former was investigated through test results on four different topics that were considered the most challenging ones applied to both the control and experimental groups. A validated, multiple-choice type, teacher-made pre-test and post-test on these particular topics were administered to gauge and compare their academic achievements in each topic. This validation was conducted by other instructors who possess expertise in the field, have taught the subject for at least five years, and obtained distinguished results on students' feedback on teaching evaluation. 

The levels of difficulty in both the pre-test and post-test were parallel and identical. Each test consists of 25 items for each topic and the students must complete it within 120 minutes. The pre-test was administered before the start of the experiment, whereas the post-test was conducted during the final exam period. All scores were converted into percentages. To determine if any statistically significant difference exists between the pre-test mean scores,  pre-test and post-test mean scores, and the post-test mean scores of the two respondent groups, we conducted a sequence of Student's $t$-tests for dependent and independent samples. To assess learners' non-cognitive gains after they experienced learning the module using the FCP, we requested them to write their opinion on the FCP in student journals as well as to respond to an FCP perception inventory. 

Table~\ref{activity} shows the classroom activities for both groups during the pre-test period. Both groups experience both the FCP and the CTM but with different topics. The first group covers the first two topics in the CTM and the other two topics in the FCP, while the second group is the opposite, first the FCP, and then the CTM. For the post-test activity, the students dedicate the remaining time to journal writing and the completion of the FCP perception inventory.
\begin{table}[H]
\tablesize{\footnotesize}
\caption{Group activities and their associated topics in college algebra conducted during the pre-test period. The boldface letter {\bfseries G} refers to the group, which can be 1 or 2, whereas {\bfseries F} and {\bfseries T} refer to the classroom approaches of flipped and traditional, respectively.} \label{activity}
\newcolumntype{C}{>{\centering\arraybackslash}X}
\begin{tabularx}{\textwidth}{CCC}
\toprule	
\textbf{Topic and Content} & \multicolumn{2}{c}{\textbf{Group} ({\bfseries G1} or {\bfseries G2})\textbf{ and Classroom Approach} ({\bfseries F} or {\bfseries T})}  \\ 
\midrule
 	    		  & {\bfseries G1T} & {\bfseries G2F} \\ 
 	    		  \cmidrule{2-3}
{\bfseries Factorization:} & \multirow{9}{4.2cm}{Students listened to the teacher's regular classroom lecture-discussion aided by some PowerPoint presentations, printed references, solved sample problems, and applications. Homework/assignment was given to the students for submission during the next meeting.} & \multirow{9}{4.2cm}{Students viewed downloaded and/or teacher-provided video and PowerPoint presentations, read modules and textbooks, and completed worksheets and other assessment materials prior to class. \linebreak Interactive activities were done by the students in the classroom.} \\
Common Monomial Factor, &    &   \\
Factors of General Trinomials, & & \\
Difference of Two Squares, & & \\
Perfect Square Trinomial, & & \\
Sum/Difference of Two Cubes & & \\ 
\cmidrule{1-1}
{\bfseries Rational Expressions:} & & \\
Operations on Rational \linebreak Expressions, & & \\
Simplifying Complex Rational Expressions & & \\ 
\midrule 
 					& {\bfseries G1F}  & {\bfseries G2T} \\
 					 \cmidrule{2-3}
{\bfseries Operations on Radicals:} & \multirow{7}{4.2cm}{Students viewed, downloaded, and/or teacher-provided video and PowerPoint presentations, read  modules and textbooks,  and completed worksheets and other assessment materials prior to class. \linebreak Interactive activities were done by the students in the classroom.} & \multirow{7}{4.3cm}{Students listened to the teacher's regular classroom lecture-discussion aided by some PowerPoint presentations, printed references, solved sample problems, and applications. \linebreak Homework/assignment was given to the students for submission during the next meeting.} \\
Simplifying Radicals, & & \\
Rationalizing Denominators & & \\ 
\cmidrule{1-1}
{\bfseries Solving Applied Problems:} & & \\
Routine and Nonroutine \linebreak 
Problems (Number, Age, & & \\
Investment, and Mensuration) \vspace{16pt}& \vspace{16pt}&\vspace{16pt} \\
\bottomrule
\end{tabularx}
\end{table}

\section{Result and Discussion}		\label{result}		\vspace*{2mm}

\subsection{Cognitive Gains}		\vspace*{1mm}

To achieve the objectives of the study, the pre-test and post-test scores were subjected to statistical treatment. Moreover, the students' journal report and their responses to the perceptions inventory were analyzed. Below are the tables showing the descriptive statistics for the gathered data and the corresponding discussion.

As reflected in Table~\ref{pre-test}, the pre-test mean scores of Group 1 and Group 2 are far below the passing score of 60\% in all topics under consideration. Very clearly, the students had very low mean scores in the topics of radical operations and solving applied problems. The small values of the respective standard deviations indicate little variation among the scores; that is, the groups are more or less homogeneous. The non-statistically significant difference between the computed $t$-values of the pre-test means confirmed the comparability of the two groups at the beginning of the experiment in so far as performance in all the topics under consideration is concerned. This indicates that any changes in their achievement can be attributed to the utilized teaching approach.

\begin{table}[H]
\caption{Descriptive statistics comparing the pre-test mean scores of the two groups (df $= 53$, $t_\textmd{\footnotesize crit} = 2.0057$). All results were not statistically significant $(p > 0.05)$.}	\label{pre-test}
\begin{tabularx}{\textwidth}{m{3cm}<{\centering}CCCCCC}
\toprule	
\multirow{2}{*}{\textbf{Topic}} & \multicolumn{2}{c}{\textbf{Group 1} \boldmath ($n = 24$)}  & \multicolumn{2}{c}{\textbf{Group 2} \boldmath ($n = 31$)} &  \multicolumn{2}{c}{\textbf{Statistics}}  \\ 
\cmidrule {2-7}
      				  & \boldmath $\bar{x}$ & \boldmath $s$ 	 &\boldmath $\bar{x}$ & \boldmath $s$ & \boldmath $|t|$-\textbf{Value} &\boldmath $p$-\textbf{Value} \\ 
      				  \midrule 
Factorization 		  & 45.50 & 11.77 	&  45.87 & {\;} 7.67 &  0.1408  & 0.8886 \\      
Rational expressions &  44.71 & {\;\:}9.72 &  45.52 & {\;} 6.27 &  0.3746  & 0.7095 \\      
Radical operations	 &  33.21 & 10.18 	&  34.35 & 10.17   &  0.4121  & 0.6819 \\      
Applied problems	 &  38.96 & {\;\:}9.29 &  39.23 & 10.46   &  0.0996  & 0.9210 \\ 
\bottomrule
\end{tabularx}
\end{table}

The pre-test and post-test mean values of each group in all topics as shown in Table~\ref{prepost-test} indicate that learning took place using any of the two approaches---flipped or traditional. In teaching factorization and rational expressions, the use of the CTM for Group 1 accounted for respective mean gains of 39.58 and 27.58 as compared with the respective mean gains of 27.52 and 36.16 which can be attributed to the use of the FCP for Group 2. Moreover, the use of the CTM in teaching radical operations and solving applied problems accounted for mean gains of 28.68 and 30.21, respectively, as compared to the respective mean gains of  32.62 and 36.42 of the students exposed to the FCP. Obviously, based on the mean gains, using the FCP is better than the CTM when it comes to the topics of rational expressions, radical operations, and solving applied problems. The highest mean gain is evident in the use of the CTM in factorization (39.58), followed by the use of the FCP in solving applied problems (36.42), rational expressions (36.16), and radical operations (32.62). Thus, the FCP proved to be a better teaching approach in three out of the four identified difficult topics in this study in terms of academic achievement gains. Statistically significant cognitive gains are also evident as indicated by the calculated $t$-values, which are all greater than the tabular $t$-values. This indicates that each approach has a significant positive effect on  learning, but the use of the FCP generally accounts for greater mean gain in the identified difficult topics in college algebra.

As observed in Table~\ref{post-test}, the result of the $t$-test for independent samples applied to the post-test mean scores indicate a statistically significant difference between the FCP and CTM in the delivery of factorization and rational expressions, but no significant difference in radical operations and solving applied problems. This implies that the CTM appears to be a better teaching approach for factorization while the FCP might be a better delivery method for rational expressions. 

\begin{table}[H]
\caption{Descriptive statistics between the pre-test and post-test mean scores of the respondents. Group~1 adopted the FCP for the third and fourth topics, with df $= 23$ and $t_\textmd{\footnotesize crit} = 2.0687$. Group~2 adopted the FCP for the first and second topics, with df $= 30$ and $t_\textmd{\footnotesize crit} = 2.0423$. All results were statistically very significant $(p < 0.001)$.}	\label{prepost-test}
{\footnotesize
\begin{adjustwidth}{-\extralength}{0cm}
\newcolumntype{C}{>{\centering\arraybackslash}X}
\begin{tabularx}{\fulllength}{m{3cm}<{\centering}CCCCCCCCCCCC}
\toprule	
\multirow{3}{*}{\textbf{Topic}}			& \multicolumn{6}{c}{\textbf{Group 1} \boldmath ($n = 24$)} & \multicolumn{6}{c}{\textbf{Group 2} \boldmath ($n = 31$)}   \\ 
\cmidrule{2-13}
					& \multicolumn{2}{c}{\textbf{Pre-Test}} & \multicolumn{2}{c}{\textbf{Post-Test}} & \multicolumn{2}{c}{\textbf{Statistics}} & \multicolumn{2}{c}{\textbf{Pre-Test}} & \multicolumn{2}{c}{\textbf{Post-Test}} & \multicolumn{2}{c}{\textbf{Statistics}} \\ 
\cmidrule{2-13} 
					&\boldmath $\bar{x}$ &\boldmath $s$ 	&\boldmath $\bar{x}$ &\boldmath $s$   & \boldmath$|t|$ &\boldmath $p$ &\boldmath $\bar{x}$ &\boldmath $s$ &\boldmath $\bar{x}$ &\boldmath $s$ & \boldmath$|t|$ &\boldmath $p$ \\
					 \midrule 
Factorization 		& 45.50 	& 11.77 & 85.08	    & 11.72 & {\!\!}3878.03   & 0.000     & 45.87     & {\;\:}7.67 & 73.39 & 10.42 & 55.72  & 0.000 \\      
Rational expressions& 44.71 	& {\;\:}9.72 &72.29     & 12.28 & {\;\:}52.78     & 0.000		& 45.52 & {\;\:}6.27 & 81.68 & 11.92 & 35.63  & 0.000 \\      
Radical operations	& 33.21 	& 10.18 & 65.83     & 10.58 & 399.51    & 0.000     & 34.35     & 10.17   & 63.03 & 14.65 & 35.64  & 0.000 \\      
Applied problems	& 38.96 	& {\;\:}9.29 &75.38 	& 14.86 & {\;\:}32.03     & 0.000     & 39.23     & 10.46   & 69.44 & 14.12 & 45.96  & 0.000 \\ 
\bottomrule
\end{tabularx}
	\end{adjustwidth}
}
\end{table}

\vspace{-6pt}

\begin{table}[H]
\caption{Descriptive statistics comparing the post-test mean scores of the two groups (df $= 53$, $t_\textmd{\footnotesize crit} = 2.0057$). Only the mean scores for the first two topics appear to be statistically significant and the FCP produced higher mean scores except for the first topic on factorization.}	\label{post-test}
\newcolumntype{C}{>{\centering\arraybackslash}X}
\begin{tabularx}{\textwidth}{m{3cm}<{\centering}CCCCCC}
\toprule	
\multirow{2}{*}{\textbf{Topic}} &  \multicolumn{2}{c}{\textbf{Group 1} \boldmath($n = 24$)}  & \multicolumn{2}{c}{\textbf{Group 2} \boldmath($n = 31$)} & \multicolumn{2}{c}{\textbf{Statistics}}  \\ 
\cmidrule{2-7}
					 & \boldmath $\bar{x}$ & \boldmath $s$ 		& \boldmath $\bar{x}$ & \boldmath $s$ 	 & \boldmath $|t|$-\textbf{Value} & \boldmath $p$-\textbf{Value} \\ 
					 \midrule 
Factorization 		 & {\!\!}85.08 	 & 11.72 	 & 73.39$^\star$ 	 & 10.42 &  3.9076  	 & {\,}0.0000$^\ast$ \\      
Rational expressions & {\!\!}72.29 	 & 12.28 	 & 81.68$^\star$ 	 & 11.92 &  2.8595  	 & {\,}0.0061$^\ast$ \\      
Radical operations	  & 65.83$^\star$ 	 & 10.58 	& {\!\!}63.03 	 & 14.65 &  0.7897  	 & {\!\!}0.4332 \\      
Applied problems	 &  75.38$^\star$ 	 & 14.86    & {\!\!}69.44 	 & 14.12 &  1.5123	     & {\!\!}0.1364 \\ 
\bottomrule
\end{tabularx}
{\footnotesize{$^\star$post-test mean score under the FCP; \quad {$^\ast$}significant at $\alpha = 0.05$ level}}
\end{table}
\vspace{-6pt}

It can be noted that the students registered a higher mean score in factorization under the CTM while the use of the FCP accounts for the higher mean scores in rational expressions, radical operations, and solving applied problems. However, non-statistically significant mean differences in the respondents' post-test scores in the last two topics imply that learning can be achieved whether delivery of these topics was done using the FCP (where they view the lessons online, through other modes, and do their supposed homework in the classroom) or the CTM (where they have their lessons in the classroom and really do their homework at home).

There could be several factors that might influence this result. The first one is related to the topic and content investigated in this study. The second factor relates to the quality of lectures in both pedagogies. The third factor is associated with the variation in students' academic performance when they encounter various college algebra topics in different learning environments, i.e., the FCP vs. CTM. Our finding suggests that learners achieved better cognitive gains on the first topic of factorization when it was delivered using the CTM instead of the FCP. This topic may be relatively less challenging than the other three topics, and for some of the contents, the students may have seen them in primary or secondary mathematics. For example, finding the greatest common factor between monomials relates to basic number theory on prime numbers in primary school mathematics. Another example is factoring general trinomials, which requires a direct explanation of how to find two integers whose product is one monomial and whose sum is another monomial for easier understanding. As might be the case in the latter, although viewing video recordings might explain the procedure, an absence of immediate feedback might hinder further progress and comprehension, and thus the CTM might be better suited for these concepts than the FCP.
	
The second factor might be related to the quality of the lectures in both traditional and flipped settings. From the learners' perspective, it does not reserve a possibility that the teaching delivery during direct instruction has better quality and the recorded video recording for this particular topic of factorization. As argued by Krantz (2015), teaching mathematics by lecturing directly---offline mode, face-to-face---and hence predominantly the CTM, is still a powerful teaching device that has stood the test of time for more than three millennia provided that the instructor does it very well~\cite{krantz2015how}. However, we need to be cautious in directly swallowing this argument since Krantz's argument was also seriously disputed by a meta-analysis study from several STEM disciplines by Freeman et al. (2014) who discovered that student academic achievements are significantly better when some kind of active learning is implemented~\cite{freeman2014active}.

For the third factor, we observed that the pre-test and post-test results from Group 1 on the first topic factorization yield sample standard deviation values of around 11.75, i.e., 11.77 and 11.72 for the former and latter, respectively. These not-so-low values suggest that the CTM only improved the students' scores and translated the mean to a higher score but it did not really close the difference gap between the more academically and less academically prepared students. In other words, the CTM has successfully prevented the widening gap between the various spectrum of students' academic strengths. By looking at the results for Group 2, the standard deviation values increase from 7.67 to 10.42 for the pre-test and post-test results, respectively. This nearly 3\% increase in standard deviation could influence a lower post-test mean score compared to Group 1, and thus the mean difference for post-test results appears to be statistically significant. This could suggest that the FCP fails to prevent gap widening or enclosing the existing gap in academic achievement for this particular topic.

A similar outcome was observed for the third topic of radical operations, but then the pedagogy was reversed. Although the mean difference for the post-test scores did not appear to be statistically significant, we could argue that the FCP has had relative success in preventing the widening gap in academic performance among students in Group 1 for this particular topic. The situation, however, was slightly different for the fourth topic of applied problems. The increases in sample standard deviations from the pre-test to post-test results were nearly 6\% and 4\% for Groups 1 and 2, respectively. This finding might suggest that the mean difference for the post-test scores did not appear to be statistically significant for this topic even though the result from the FCP in Group 1 was better than the one from the CTM in Group 2. A similar argument could be applied to the second topic of rational expressions where both groups experienced an increase in standard deviation values. But since the mean score for Group 2 is much higher than Group 1 (nearly 9\%), the increase in standard deviations was still compensated and the mean difference for the post-test scores appeared to be statistically significant. It would be interesting to investigate whether other studies yield comparable outcomes when a similar methodology was applied to different cohorts of students.

\subsection{Non-Cognitive Gains: Attitudes Toward Mathematics}		\vspace*{1mm}

In addition to cognitive gains of better academic performance, the students also acquired non-cognitive gains after experiencing the FCP. By soliciting their opinions regarding the use of the FCP through a flipped classroom perception inventory, we collect further insight regarding their learning characteristics. Table~\ref{perception} summarizes the percentages of those who agreed on each item in the perception inventory.

From this outcome, we observe that students' responses are generally favorable when it comes to learning using the FCP. All the ten items yield either more than 90\% or nearly 90\% of students' ratings, where the highest percentage (94\%) belongs to the seventh item that the FCP promotes a positive attitude toward mathematics. The high percentages of students in favor of each item in the inventory indicate that they welcome this type of pedagogy in learning mathematics, particularly college algebra. In addition to improving positive attitudes toward mathematics, providing the opportunity to understand mathematical concepts in a better way, and offering various alternatives in dealing with mathematical problems, the non-cognitive benefit of the FCP also goes beyond mathematics learning per~se.

\begin{table}[H]
\caption{The percentage of students who agreed on each item in the FCP perception inventory. The rating is either more than 90\% or very close to 90\%.}	\label{perception}
\newcolumntype{C}{>{\centering\arraybackslash}X}
		\begin{tabularx}{\textwidth}{Cm{9cm}<{\raggedright}C}
\toprule	
\textbf{No}. & \textbf{Item} & \textbf{Rating (\%)} \\ 
\midrule
1   & Motivates students to study & 92 \\ 
2   & Develops students' confidence to solve problems & 92 \\
3   & Offers a variety of alternatives in understanding the lesson & 93 \\ 
4   & Improves students' creative and critical thinking ability & 92 \\
5   & Provides students the opportunity to understand mathematical concepts better & 92 \\ 
6   & Allows students to participate actively in the learning process and progress independently & 91 \\
7   & Promotes positive attitudes toward mathematics & 94 \\ 
8   & Strengthens students' retention of subject matter & 88 \\
9   & Promotes better cooperation between teacher and students and among students & 92 \\ 
10  & Offers powerful ways of dealing with mathematical problems & 89 \\
\bottomrule
\end{tabularx}
\end{table}

In the essence of SLE, the students in our study confirmed that the FCP has motivated them to study, improved their critical thinking, enhanced their knowledge retention rate, and promoted better communication with the teachers and among themselves. Overall, the FCP has assisted them not only in acquiring non-cognitive gains but also in appreciating some aspects of SLE. Note that although some of them may forget the details of certain mathematical concepts or theorems in this particular course, the acquired qualities and positive attitudes will be transferable to another advanced or different course, even to a totally different setting, whether in their profession as educators or elsewhere; see~\cite{ben2021sustainable}.

Other non-cognitive gains from the use of the FCP indicated by the students in their journals are as follows: a higher level of participation in active learning sessions, a stronger sense of being connected to the teacher and with other students, a positive perception of the course, and an enjoyable learning environment. Moreover, they acquired from wide resources of learning tools at varying speeds in their chosen learning environment and were able to have a continuing review as the need arises in as much as the materials were accessible even for absentees. Hence, the element of flexibility. Finally, they also like the idea that their teacher had more time in guiding them during activity sessions.

Empirical results confirmed that the use of the FCP enhances students' intellectual experiences and learning outcomes. Additionally, although some students in the CTM setting were more satisfied with the clarity of instruction provided by the teacher and the facility of getting feedback immediately after the lecture, they felt more strongly in terms of gaining a greater appreciation of college algebra concepts through the FCP.

Some unedited comments culled from the students' journals are given as follows:\\

\begin{quote}
\emph{ ``Learning math has never been like this before because this time I did not have to sit in class very long doing nothing but listen and be bored.'';}
\end{quote}

\begin{quote}
\emph{``I enjoyed learning math because I can switch from viewing videos to reading the PowerPoint presentation sent by my teacher to simply reading the textbook, trying to get ready for our class activities.'';}
\end{quote}

\begin{quote}
\emph{``Whenever I have questions, I kept on viewing the online lesson wherever I am so long as I have my laptop with me and I have access to internet.'';}
\end{quote}

\begin{quote}
\emph{``While viewing/reading the lesson, I had to list down my questions and  ask my teacher and my classmates for clarification.'';}
\end{quote}

\begin{quote}
\emph{``I could have learned more math had this approach been introduced earlier.'';}
\end{quote}

\begin{quote}
\emph{``I felt helpless whenever I had questions, but nobody was there to answer immediately.'';}
\end{quote}

\begin{quote}
\emph{``Sometimes, I cannot access the lesson because there was no internet facility.'';}
\end{quote}

\begin{quote}
\emph{``There should be a larger viewing center at the library where we can stay whenever no internet is available at home.'';}
\end{quote}

\begin{quote}
\emph{``Sometimes I was tempted to skip my class when I thought I already understand the lesson very well and my answers to the assessment materials were almost alright.'';}
\end{quote}

\begin{quote}
\emph{``I enjoyed the extra time given me to do the viewing of the lesson over and over again until it becomes clear to me.'';}
\end{quote}

\begin{quote}
\emph{``I got a strong bonding session with my classmates and my teacher during the face-to-face encounter.'';}
\end{quote}

\begin{quote}
\emph{``As I read my lessons alone, I  became independent, patient and confident in learning math.'';}
\end{quote}

\begin{quote}
\emph{``It seems that I became more interested in math.'';}
\end{quote}

\begin{quote}
\emph{``During class time I had a fruitful exchange of ideas with my peers as well as with our~teacher.''}
\end{quote}

From these comments, we observe that there is a mixed reaction from the students regarding their experience with the FCP. For the negative aspect of the situation, we could categorize it into at least three factions: the students' side, the teachers' part, and the intermediary party that links the learners with their facilitators. For the first group, some students consider that studying on their own is a mammoth challenge while for others, they could understand easily. The former might feel helpless whenever they have some questions that need immediate assistance, but nobody was there to provide prompt help and feedback. For the latter, they might be tempted to skip the face-to-face sessions with their teachers. For both sides of the spectrum, we observe that the students might be at risk, and thus educators should interfere promptly and determine better strategies to facilitate their sustainable learning.

The challenge in the intermediary party that links students with their teachers is mostly related to technical issues, such as Internet connectivity and computer facilities. Several students brought up these points in the aforementioned comments, and their learning might be hindered as well. This leads us to address burdens on teachers' shoulders. Although they are not mentioned in the students' comments, the previous two challenges would commit educators to additional tasks in assisting their learners. They might need to prepare additional teaching materials for the students who were lacking with Internet access, they might need to schedule additional activities for the students who struggle understanding basic materials, or they might adjust the teaching speed and material coverage to accommodate those who were left behind.

Despite these downsides, we recognize that more students wrote about their constructive experience with the FCP. They enjoyed learning mathematics and became more interested in delving deeper into mathematics. In addition to acting more independently when studying on their own, they also developed patience and gained more confidence in learning mathematics. Since they underwent different learning experiences, they did not feel bored anymore when studying mathematics. Clearly, many of these students have an improved, better attitude toward mathematics; we hope that their learning will last since some of them will eventually become educators themselves. This type of empowering non-cognitive gains is transferable not only to this particular module but also to other skills they might learn in the future. To a great extent, by actively renewing and relearning through self-assessment and inquiry, they are on the right track to adopting SLE; see~\cite{ben2021sustainable}.
 
\section{Conclusions and Recommendation}		\label{conclude}		\vspace*{1mm}

We admit several limitations in this study. In addition to a relatively small number of participants, we did not collect any details such as their precise age and/or birthday as well as gender information. Another concern relates to the outcome of the experiment. Although we initially hypothesized that the FCP would always yield better learners' cognitive gains for all the four topics in college algebra considered in this study, it turned out that the first topic on factorization yielded better students' cognitive gains under the CTM. A further investigation on the root of this cause might be an interesting topic and a potentially open problem that is worth pursuing further, as we recommend in the final paragraph of this~section.

To conclude, we have considered an attempt in improving essential competencies for sustainability, particularly in the context of mathematics teaching and learning in higher education. These ESD competencies are not only crucial but also cross-cut for achieving other goals in SDGs in addition to attaining SDG4. The focus of our study is cultivating sustainable learning among learners in mathematics education through active learning and the student-centered FCP. Learning formats under this pedagogy allow learners not only to renew and rebuild knowledge but also to cope with challenging and complex situations. In particular, we have selected students majoring in mathematics education who were enrolled in a college algebra module at a private university in the Philippines. Since the majority of these students will become educators themselves, they need to acquire the habit of lifelong learning early in their careers, thus adopting sustainable learning for themselves and cultivating it to their potential pupils. 
	
Furthermore, since four topics were considered pretty challenging for many students enrolled in this course, our quantitative study then concentrated on the assessment results of these four topics. For both the control (CTM) and experimental (FCP) groups, the cognitive gains across the four topics in college algebra appeared to be statistically significant, as indicated by improved mean scores in the post-test. However, by comparing data of the post-test results, we observed that the experimental Group 1 achieved better cognitive gains under the FCP for the third and fourth topics of radical operations and applied problems, respectively, although the mean difference with the control Group 2 under the CTM did not appear to be statistically significant. The outcome for the second topic materialized as we have anticipated, Group 2 under the FCP gained better cognitively than Group 1 under CTM. However, for the first topic, the outcome is rather surprising where cognitive gains under the CTM trumped those of the FCP. In both instances, the mean differences appeared to be statistically significant. This seemingly rather puzzling outcome could potentially invite further discussion, debate, and other follow-up studies related to this topic.

Regarding non-cognitive gains, we qualitatively analyzed the result of the FCP perception inventory as well as student journals that they wrote regarding their feedback on the FCP. For the former, the rating for each of the ten items is overwhelmingly positive, ranging from 88\% to 94\%. In particular, promoting a positive attitude toward mathematics gained the highest rate of 94\%. For the latter, we observed that the reaction from the students is mixed although they also generally wrote positive things about the FCP. Some of the hindrances are related to technical issues, such as Internet connectivity. However, many students were happier to have control of their own learning, more confident, learned more, enjoyed learning, and developed more interest in mathematics. Certainly, these positive traits could lead to more active and intentional learning, not only in mathematics modules but also in other courses. Any incremental skills and strategies that they acquired in this course could also potentially be transferred beyond the classroom setting, into different contexts and domains. Hence, these qualitative findings suggest that the FCP is a possible excellent teaching approach for cultivating SLE among mathematics learners. 

Upon conducting this study, on the one hand, we would not recommend entirely abandoning the CTM and its predominantly lecturing method in teaching mathematics at all levels of education. Lectures and the CTM, after all, have their place in transmitting and imparting mathematical knowledge to the learners albeit their nature tends to be teacher-centered and students' role is passive-receptive. Coupled with the fact that many teachers do not lecture very well when they teach, or rather, they do not teach very well when they lecture, it is imperative to equip mathematics educators---whether they are teachers at the primary and secondary levels, or instructors and professors at the tertiary level---with sufficient training in teaching techniques, including delivering lectures. Marrying excellent lecturing with active learning approaches, such as the FCP, may have a tremendous positive impact on students' learning. When learners' curiosity is aroused and they experience a sparkle of light in mathematical understanding, they will take the next step of actively learning beyond the course syllabus, renewing their learning, whether independently or collaboratively, and transferring it to different domains. Indeed, they are on the right track toward lifelong learning and SLE.

On the other hand, we would like to recommend practitioners in education to continue embracing and implementing student-centered, active learning approaches in their classrooms, such as the FCP. Although there are some cases where the FCP fails in enhancing students' learning and arousing their interest in learning, the majority of studies in the body of published literature consistently support positive outcomes of the FCP, along with other active and sustainable learning methodologies, such as problem-based and/or project-based learning. Very often, when implementing the FCP, some instructors were distracted from the main objectives of the pedagogy itself by focusing too much on preparing video recordings that the students are supposed to view before they come to class and hence they come unprepared to the class. However, we need to emphasize that the main point of the FCP is what happens when learners are inside the class, doing and solving problems, collaborating with their peers, and interacting with the instructor. When learners take charge and are actively involved in their own learning, they are one step closer to embracing sustainable learning for their own life.

To end this recommendation, we would like to invite other researchers to delve deeper into some possible causes regarding the effectiveness of the FCP vs. the CTM when it comes to cognitive gains. Since our results suggest different outcomes for different algebraic topics, we would like to see a similar study be replicated in the future, preferably on a wider scale of subjects and objects. Any potential outcomes would be interesting for further discussion, debate, and collaboration.

\vspace{6pt}
\authorcontributions{Conceptualization, M.J.A.; methodology, M.J.A.; software, M.J.A. and N.K.; validation, N.K. and M.J.A.; formal analysis, M.J.A. and N.K.; investigation, M.J.A.; resources, M.J.A.; data curation, M.J.A.; writing---original draft preparation, M.J.A.; writing---review and editing, N.K.; visualization, N.K.; supervision, M.J.A.; project administration, M.J.A.; funding acquisition, N.K. All authors have read and agreed to the published version of the manuscript.}

\funding{This research was supported by the National Research Foundation (NRF) of Korea and funded by the Korean Ministry of Science, Information, Communications, and Technology (MSICT) through Grant No.~NRF-2022-R1F1A-059817 under the scheme of Broadening Opportunities Grants---General Research Program in Basic Science and~Engineering.}

\institutionalreview{Not applicable.}

\informedconsent{Not applicable.}

\dataavailability{Data available upon request.} 

\acknowledgments{The authors acknowledged all three anonymous reviewers whose comments and remarks have improved the article significantly.}

\conflictsofinterest{The authors declare no conflict of interest.}


\begin{adjustwidth}{-\extralength}{0cm}
\reftitle{References}
	
\end{adjustwidth}	

\begin{thebibliography}{99}
\bibitem{un2030sdg} Department of Economic and Social Affairs, United Nations. \ Sustainable Development Goals. 2015. Available online: \url{https://www.un.org/development/desa/disabilities/envision2030.html} (accessed on 30 July 2022). 

\bibitem{hesi2021higher} Department of Economic and Social Affairs, UN. The HESI Working Group. \ Higher Education Sustainability Initiative, Assessments of Higher Education's Progress towards the UN Sustainable Development Goals, Volume 2: For Higher Education Institutions Participating in Assessments. 2021. Available online: \url{https://sdgs.un.org/} (accessed on 30 July 2022).

\bibitem{unesco2007education} UNESCO. \ Education for Sustainable Development Goals: Learning Objectives. 2017. Available online: \url{https://unesdoc.unesco.org/ark:/48223/pf0000247444} (accessed on 30 July 2022).

\bibitem{rogovaya2019critical} Rogovaya, O.S.; Larchenkova, L.; Gavronskaya, Y. \ Critical thinking in STEM (science, technology, engineering, and mathematics). Utopía y Praxis Latinoamericana: Revista Internacional de Filosofía Iberoamericana y Teoría Social. \emph{Lat. Am. Utop. Praxis: Int. J. -Ibero-Am. Philos. Soc. Theory} \textbf{2019}, \emph{24}, 32--41. [\href{https://dialnet.unirioja.es/servlet/articulo?codigo=7406882}{Dialnet}]

\bibitem{kioupi2019education} Kioupi, V.; Voulvoulis, N. \ Education for sustainable development: A systemic framework for connecting the SDGs to educational outcomes. \emph{Sustainability} \textbf{2019}, \emph{11}, 6104. [\href{http://dx.doi.org/10.3390/su11216104}{CrossRef}] 

\bibitem{galan2020sustainable} Galán-Casado, D.; Moraleda, A.; Martínez-Martí, M.L.; Pérez-Nieto, M.Á. \ Sustainable environments in education: Results on the effects of the new environments in learning processes of university students. \emph{Sustainability} \textbf{2020}, \emph{12}, 2668. [\href{http://dx.doi.org/10.3390/su12072668}{CrossRef}]

\bibitem{collado2021flipped} Collado-Valero, J.; Rodríguez-Infante, G.; Romero-González, M.; Gamboa-Ternero, S.; Navarro-Soria, I.; Lavigne-Cerván, R. \ Flipped classroom: Active methodology for sustainable learning in higher education during social distancing due to COVID-19. \textit{Sustainability} \textbf{2021}, \emph{13}, 5336. [\href{http://dx.doi.org/10.3390/su13105336}{CrossRef}]

\bibitem{ben2021sustainable} Ben-Eliyahu, A. \ Sustainable learning in education. \emph{Sustainability} \textbf{2021}, \emph{13}, 4250. [\href{http://dx.doi.org/10.3390/su13084250}{CrossRef}]

\bibitem{venkataraman2009education} Venkataraman, B. \ Education for sustainable development. \emph{Environ. Sci. Policy Sustain. Dev.} \textbf{2009}, \emph{51}, 8--10. [\href{http://dx.doi.org/10.3200/ENVT.51.2.08-10}{CrossRef}]

\bibitem{hopkins2002education} Hopkins, C.; McKeown, R. \ Education for sustainable development: An international perspective. In \emph{Education and Sustainability: Responding to the Global Challenge}; Tilbury, D., Stevenson, R.B., Fien, J., Schreuder, D., Eds.; \ International Union for Conservation for Nature and Natural Resources--Commission on Education and Communication: Gland, Switzerland; Cambdridge, UK, 2002; pp.~13--24. 

\bibitem{boeve2015effectiveness} Boeve-de, Pauw, J.; Gericke, N.; Olsson, D.; Berglund, T. \ The effectiveness of education for sustainable development. \emph{Sustainability} \textbf{2015}, \emph{7}, 15693--15717. [\href{https://doi.org/10.3390/su71115693}{CrossRef}]

\bibitem{huckle2015UN} Huckle, J.; Wals, A.E. \ The UN decade of education for sustainable development: Business as usual in the end. \emph{Environ. Educ. Res. }\textbf{2015}, \emph{21}, 491--505. [\href{http://dx.doi.org/10.1080/13504622.2015.1011084}{CrossRef}]

\bibitem{graham2015sustainable} Graham, L.; Berman, J.; Bellert, A. \ \emph{Sustainable Learning}; Cambridge University Press: Cambridge, UK, 2015.

\bibitem{hays2020sustainable} Hays, J.; Reinders, H. \ Sustainable learning and education: A curriculum for the future. \emph{Int. Rev. Educ.} \textbf{2020}, \emph{66}, 29--52. [\href{http://dx.doi.org/10.1007/s11159-020-09820-7}{CrossRef}]

\bibitem{moore2005seven} Moore, J. \ Seven recommendations for creating sustainability education at the university level: A guide for change agents. \emph{Int. J. Sustain. High. Educ.} \textbf{2005}, \emph{6}, 326--339. [\href{http://dx.doi.org/10.1108/14676370510623829}{CrossRef}]

\bibitem{zimmerman1990self} Zimmerman, B.J. \ Self-regulated learning and academic achievement: An overview. \emph{Educ. Psychol.} \textbf{1990}, \emph{25}, 3--17. [\href{http://dx.doi.org/10.1207/s15326985ep2501_2}{CrossRef}]

\bibitem{pintrich1995understanding} Pintrich, P.R. \ Understanding self‐regulated learning. \emph{New Dir. Teach. Learn.} \textbf{1995}, \emph{1995}, 3--12. [\href{http://dx.doi.org/10.1002/tl.37219956304}{CrossRef}]

\bibitem{boekaerts1999self} Boekaerts, M. \ Self-regulated learning: Where we are today. \emph{Int. J. Educ. Res}. \textbf{1999}, \emph{31}, 445--457. [\href{http://dx.doi.org/10.1016/S0883-0355(99)00014-2}{CrossRef}]

\bibitem{puustinen2001models} Puustinen, M.; Pulkkinen, L. \ Models of self-regulated learning: A review. \emph{Scand. J. Educ. Res.} \textbf{2001}, \emph{45}, 269--286. [\href{http://dx.doi.org/10.1080/00313830120074206}{CrossRef}]

\bibitem{zimmerman2011self} Zimmerman, B.J.; Schunk, D.H. \ Self-regulated learning and performance: An introduction and an overview. In \emph{Handbook of Self-Regulation of Learning and Performance}; Schunk, D.H., Zimmerman, B., Eds.; Routledge: New York, NY, USA, 2011; pp.~15--26. 

\bibitem{cortese2003critical} Cortese, A.D. \ The critical role of higher education in creating a sustainable future. \emph{Plan. High. Educ.} \textbf{2003}, \emph{31}, 15--22. Available online: \url{https://www.redcampussustentable.cl/wp-content/uploads/2022/07/6-CorteseCriticalRoleOfHE.pdf} (accessed 30 July 2022).

\bibitem{jones2010sustainability} Jones, P.; Selby, D.; Sterling, S. (Eds.) \ \emph{Sustainability Education: Perspectives and Practice Across Higher Education}; Routledge: Oxfordshire, UK, 2010.

\bibitem{zaleniene2021higher} Žalėnienė, I.; Pereira, P. \ Higher education for sustainability: A global perspective. \emph{Geogr. Sustain.} \textbf{2021}, \emph{2}, 99--106. [\href{http://dx.doi.org/10.1016/j.geosus.2021.05.001}{CrossRef}]

\bibitem{veidemane2022education} Veidemane, A. \ Education for sustainable development in higher education rankings: Challenges and opportunities for developing internationally comparable indicators. \emph{Sustainability} \textbf{2022}, \emph{14}, 5102. [\href{http://dx.doi.org/10.3390/su14095102}{CrossRef}]

\bibitem{sterling2001sustainable} Sterling, S.; Orr, D. \ \emph{Sustainable Education: Re-visioning Learning and Change}; Green Books for the Schumacher Society: Totnes, UK, 2001.

\bibitem{branden2012sustainable} Van den Branden, K. \ Sustainable education: Basic principles and strategic recommendations. \emph{Sch. Eff. Sch. Improv.} \textbf{2012}, \emph{23}, 285--304. [\href{http://dx.doi.org/10.1080/09243453.2012.678865}{CrossRef}]

\bibitem{fullan2015new} Fullan, M. \ \emph{The New Meaning of Educational Change}, 5th ed.; Teachers College Press: New York, NY, USA, 2015.

\bibitem{peris2015sustainable} Peris-Ortiz, M.; Lindahl, J.M.M. (Eds.) \ \emph{Sustainable Learning in Higher Education: Developing Competencies for the Global Marketplace}; Springer: Berlin/Heidelberg, Germany, 2015.

\bibitem{sa2018transversal} Sá, M.J.; Serpa, S. \ Transversal competences: Their importance and learning processes by higher education students. \emph{Educ. Sci.} \textbf{2018}, \emph{8}, 126. [\href{http://dx.doi.org/10.3390/educsci8030126}{CrossRef}]

\bibitem{membrillo2021sustainability} Membrillo-Hernández, J.; Lara-Prieto, V.; Caratozzolo, P. \ Sustainability: A public policy, a concept, or a competence? Efforts on the implementation of sustainability as a transversal competence throughout higher education programs. \emph{Sustainability} \textbf{2021}, \emph{13}, 13989. [\href{http://dx.doi.org/10.3390/su132413989}{CrossRef}]

\bibitem{renert2011mathematics} Renert, M. \ Mathematics for life: Sustainable mathematics education. \emph{Learn. Math.} \textbf{2011}, \emph{31}, 20--26. [\href{https://www.jstor.org/stable/41319547}{JSTOR}]

\bibitem{edwards2010organizational} Edwards, M. \ \emph{Organizational Transformation for Sustainability: An Integral Metatheory}; Routledge: New York, NY, USA, 2010.

\bibitem{summer2020sustainable} Summer, A. \ A sustainable way of teaching basic mathematics. \emph{Discourse Commun. Sustain. Educ.} \textbf{2020}, \emph{11}, 106--120. [\href{http://dx.doi.org/10.2478/dcse-2020-0021}{CrossRef}]

\bibitem{hamilton2014sustainability} Hamilton, J.; Pfaff, T.J. \ Sustainability education: The what and how for mathematics. \emph{PRIMUS} \textbf{2014}, \emph{24}, 61--80. [\href{http://dx.doi.org/10.1080/10511970.2013.834526}{CrossRef}]

\bibitem{barwell2018some} Barwell, R. \ Some thoughts on a mathematics education for environmental sustainability. In \emph{The Philosophy of Mathematics Education Today}; ICME-13 Monographs; Ernest, P., Ed.; Springer: Cham, Switzerland, 2018; pp. 145--160.

\bibitem{widiati2019philosophy} Widiati, I.; Juandi, D. \ Philosophy of mathematics education for sustainable development. \emph{J. Physics: Conf. Ser.} \textbf{2019}, \emph{1157}, 022128. [\href{http://dx.doi.org/10.1088/1742-6596/1157/2/022128}{CrossRef}]

\bibitem{li2020philosophy} Li, H.C.; Tsai, T.L. \ Philosophy of education for sustainable development in mathematics education: Have we got one? \emph{J.~Math.-Teach.-Res. J.} \textbf{2020}, \emph{12}, 136--140. [\href{https://commons.hostos.cuny.edu/mtrj/wp-content/uploads/sites/30/2020/09/v12n2-Philosophy-of-education-for-sustainable-development.pdf}{CUNY}]

\bibitem{li2021education} Li, H.C.; Tsai, T.L. \ Education for sustainable development in mathematics education: What could it look like? \emph{Int. J. Math. Educ. Sci. Technol.} \textbf{2021}, \emph{53}, 2532--2542. [\href{http://dx.doi.org/10.1080/0020739X.2021.1941361}{CrossRef}]

\bibitem{moreno2022training} Moreno-Pino, F.M.; Jiménez-Fontana, R.; Domingo, J.M.C.; Goded, P.A. \ Training in mathematics education from a sustainability perspective: A case study of university teachers' views. \emph{Educ. Sci.} \textbf{2022}, \emph{12}, 199. [\href{http://dx.doi.org/10.3390/educsci12030199}{CrossRef}]

\bibitem{lage2000inverting} Lage, M.J.; Platt, G.J.; Treglia, M. \ Inverting the classroom: A gateway to creating an inclusive learning environment. \emph{ J. Econ. Educ.} \textbf{2000}, \emph{31}, 30--43. [\href{http://dx.doi.org/10.1080/00220480009596759}{CrossRef}]

\bibitem{connors2000thayer} Connors, E. \ The Thayer method: Student active learning with positive results. \emph{J. Math. Sci. Collab. Explor}. \textbf{2000}, \emph{4}, 101--117. [\href{https://doi.org/10.25891/P6SW-GR81}{CrossRef}]

\bibitem{shell2002thayer} Shell, A.E. \ The Thayer method of instruction at the United States Military Academy: A modest history and a modern personal account. \emph{Probl. Resour. Issues Math. Undergrad. Stud.} \textbf{2002}, \emph{12}, 27--38. [\href{http://dx.doi.org/10.1080/10511970208984015}{CrossRef}]

\bibitem{overmyer2012flipped} Overmyer, J. \ Flipped Classrooms 101. Principal (September/October), 2012. pp. 46--47. Available online: \href{https://www.naesp.org/sites/default/files/Overmyer_SO12.pdf}{\url{https://www.naesp.org/sites/default/files/Overmyer_SO12.pdf}} (accessed on 30 July 2022).

\bibitem{king1993from} King, A. \ From sage on the stage to guide on the side. \emph{Coll. Teach.} \textbf{1993}, \emph{41}, 30--35. [\href{http://dx.doi.org/10.1080/87567555.1993.9926781}{CrossRef}]

\bibitem{fox2001constructivism} Fox, R. \ Constructivism examined. \emph{Oxf. Rev. Educ}. \textbf{2001}, \emph{27}, 23--35. [\href{http://dx.doi.org/10.1080/03054980125310}{CrossRef}]

\bibitem{arends1998resource} Arends, R.I. \ \emph{Resource Handbook. Learning to Teach}, 4th ed; McGraw-Hill: Boston, MA, USA, 1998.

\bibitem{philips1995good} Phillips, D.C. \ The good, the bad, and the ugly: The many faces of constructivism. \emph{Educ. Res.} \textbf{1995}, \emph{24}, 5--12. [\href{http://dx.doi.org/10.3102/0013189X024007005}{CrossRef}]

\bibitem{good1994looking} Good, T.L.; Brophy, J.E. \ \emph{Looking in Classrooms}; HarperCollins College Publishers: New York, NY, USA, 1994.

\bibitem{cooperstein2004beyond} Cooperstein, S.E.; Kocevar‐Weidinger, E. \ Beyond active learning: A constructivist approach to learning. \emph{Ref. Serv. Rev.} \textbf{2004}, \emph{32}, 141--148. [\href{http://dx.doi.org/10.1108/00907320410537658}{CrossRef}]

\bibitem{elliott2000educational} Elliott, S.N.; Kratochwill, T.R.; Littlefield, C.J.; Travers, J. \ \emph{Educational Psychology: Effective Teaching, Effective Learnin}g, 3rd ed.; McGraw-Hill: Boston, MA, USA, 2000.

\bibitem{vigotsky1978mind} Vygotsky, L.S. \ \emph{Mind in Society: The Development of Higher Psychological Processes}; Harvard University Press: Cambridge, MA, USA, 1978.

\bibitem{strayer2007effects} Strayer, J. \ {The Effects of the Classroom Flip on the Learning Environment: A Comparison of Learning Activity in a Traditional Classroom and a Flip Classroom That Used an Intelligent Tutoring System}. Ph.D. Thesis, The Ohio State University, Columbus, OH, USA, 2007 . Available online: \url{https://etd.ohiolink.edu/apexprod/rws_etd/send_file/send?accession=osu1189523914&disposition=inline} (accessed on 30 July 2022).

\bibitem{bretzmann2013how} Bretzmann, J. \ \textit{Flipping 2.0: Practical Strategies for Flipping Your Class}; Bretzmann Group: New Berlin, WI, USA, 2013.

\bibitem{bergmann2012flip} Bergmann, J.; Sams, A. \ \textit{Flip Your Classroom: Talk To Every Student In Every Class Every Day}; International Society for Technology in Education, Eugene, OR, USA, 2012.

\bibitem{tucker2012flipped} Tucker, B. \ The flipped classroom. \emph{Educ. Next} \ \textbf{2012}, \emph{12}, 82--83. Available online: \url{http://www.msuedtechsandbox.com/MAETELy2-2015/wp-content/uploads/2015/07/the_flipped_classroom_article_2.pdf} (accessed on 30 July 2022).

\bibitem{fulton2012upside} Fulton, K. \ Upside down and inside out: Flip your classroom to improve student learning. \textit{Learn.~Lead.~Technol}. \textbf{2012}, \emph{39}, 12--17. [\href{https://eric.ed.gov/?id=EJ982840}{ERIC}]

\bibitem{herreid2012case} Herreid, C.F.; Schiller, N.A. \ Case studies and the flipped classroom. \textit{J. Coll. Sci. Teach.} \textbf{2013}, \emph{42}, 62--66. [\href{https://www.jstor.org/stable/43631584}{JSTOR}] 

\bibitem{bishop2013flipped} Bishop, J.; Verleger, M.A. \ The flipped classroom: A survey of the research. In Proceedings of the {The 120th ASEE Annual Conference \& Exposition Conference Proceedings}, Atlanta, GA, USA, 23--26 June 2013; pp. 23--39.

\bibitem{hantla2014effects} Hantla, B.F. \ {The Effects of Flipping the Classroom on Specific Aspects of Critical Thinking in a Christian College: A Quasi-Experimental, Mixed-Methods Study}. Ph.D. Thesis, Southeastern Baptist Theological Seminary, Wake Forest, NC, USA, 2014. Available online: \url{https://www.proquest.com/docview/1547356249?pq-origsite=gscholar&fromopenview=true}  (accessed on 30 July 2022). 

\bibitem{bergman2016flipped} Bergmann, J.; Sams, A. \ \emph{Flipped Learning for Elementary Instruction}; International Society for Technology in Education: Eugene, OR, USA, 2016.

\bibitem{alrouqi2019flipped} Alrouqi, F. \ Using Flipped Classrooms to Teach Mathematics to Elementary Students in Saudi Arabia. Ph.D. Thesis, University of South Florida, Tampa, FL, USA, 2019. Available online: \url{https://www.proquest.com/docview/2321034820?pq-origsite=gscholar&fromopenview=true}  (accessed on 30 July 2022). 

\bibitem{panahi2019reviewing} Panahi, M.; Jafarkhani, F.; Bozorg, Z.J.; Nikkho, L. \ Reviewing learning environments: Effect of flipped classroom on learning levels of mathematics in primary schools. In \emph{Proceedings of the ICERI2019 12th International Conference of Education, Research and Innovation, Seville, Spain, 11--13 November 2019}; International Academy of Technology, Education and Development (IATED): Valencia, Spain, 2019;  pp. 8561--8566. 

\bibitem{khalel2021employing} Khalel, I.A.; Al-tmaran, O.S.; Hashmi, A.E. \  Employing the flipped classroom strategy in primary mathematics classes. \emph{Int. J. Res. Educ. Sci.} \textbf{2021}, \emph{4}, 497--529. [\href{http://dx.doi.org/10.29009/ijres.4.1.12}{CrossRef}]

\bibitem{saunders2014flipped} Saunders, J.M. \ The Flipped Classroom: Its Effect on Student Academic Achievement and Critical Thinking Skills in High School Mathematics. Ph.D. Thesis, Liberty University, Lynchburg, VA, USA, 2014. Available online: \url{https://www.proquest.com/docview/1639087375?pq-origsite=gscholar&fromopenview=true} (accessed on 30 July 2022). 

\bibitem{bhagat2016impact} Bhagat, K.K.; Chang, C.-N.; Chang, C.-Y. \ The impact of the flipped classroom on mathematics concept learning in high school. \emph{J.~Educ. Technol. Soc.} \textbf{2016}, \emph{19}, 134--142. [\href{https://www.jstor.org/stable/pdf/jeductechsoci.19.3.134.pdf}{JSTOR}]

\bibitem{makinde2017flipped} Makinde, S.O.; Yusuf, M.O. \ The flipped classroom: Its effects on students' performance and retention in secondary school mathematics classroom. \emph{Int. J. Innov. Technol. Integr. Educ.} \textbf{2017}, \emph{1}, 117--126.

\bibitem{lo2017using} Lo, C.K.; Hew, K.F. \ Using ``first principles of instruction'' to design secondary school mathematics flipped classroom: The findings of two exploratory studies. \emph{J. Educ. Technol. Soc.} \textbf{2017},  \emph{20}, 222--236. [\href{https://www.jstor.org/stable/pdf/jeductechsoci.20.1.222.pdf}{JSTOR}]

\bibitem{lo2017critical} Lo, C.K.; Hew, K.F. \ A critical review of flipped classroom challenges in K-12 education: Possible solutions and recommendations for future research. \emph{Res. Pract. Technol. Enhanc. Learn.} \textbf{2017}, \emph{12}, 1--22. [\href{http://dx.doi.org/10.1186/s41039-016-0044-2}{CrossRef}] [\href{https://pubmed.ncbi.nlm.nih.gov/30613253/}{PubMed}]

\bibitem{weinhandl2018technology} Weinhandl, R.; Lavicza, Z.; Süss-Stepancik, E. \ Technology-enhanced flipped mathematics education in secondary schools: A synopsis of theory and practice. \emph{K-12 STEM Educ.} \textbf{2018}, \emph{4}, 377--389. [\href{https://doi.org/10.14456/k12stemed.2018.9}{CrossRef}] [\href{https://www.learntechlib.org/p/209574/}{IPST}]

\bibitem{wei2020effect} Wei, X.; Cheng, I.; Chen, N.S.; Yang, X.; Liu, Y.; Dong, Y.; Zhai, X. \ Effect of the flipped classroom on the mathematics performance of middle school students. \emph{Educ. Technol. Res. Dev}. \textbf{2020}, \emph{68}, 1461--1484. [\href{https://doi.org/10.1007/s11423-020-09752-x}{CrossRef}]

\bibitem{overmyer2014flipped} Overmyer, G.R. \ The Flipped Classroom Model for College Algebra: Effects on Student Achievement. Ph.D. Thesis, Colorado State University, Fort Collins, CO, USA, 2014. Available online: \url{https://www.proquest.com/docview/1615100148?pq-origsite=gscholar&fromopenview=true} (accessed on 30 July 2022).

\bibitem{vansickle2015adventures} Van Sickle, J. \ Adventures in flipping college algebra. \emph{PRIMUS} \textbf{2015}, \emph{25}, 600--613. [\href{https://doi.org/10.1080/10511970.2015.1031299}{CrossRef}]

\bibitem{jungic2015flipping} Jungić, V.; Kaur, H.; Mulholl, J.; Xin, C. \ On flipping the classroom in large first year calculus courses. \emph{Int. J. Math. Educ. Sci. Technol.} \textbf{2015}, \emph{46}, 508--520. [\href{https://doi.org/10.1080/0020739X.2014.990529}{CrossRef}]

\bibitem{maciejewski2016flipping} Maciejewski, W. \  Flipping the calculus classroom: An evaluative study. \emph{Teach. Math. Its Appl. Int. J. IMA} \textbf{2016}, \emph{35}, 187--201. [\href{https://dx.doi.org/10.1093/teamat/hrv019}{CrossRef}]

\bibitem{mcgivney2013flipping} McGivney-Burelle, J.; Xue, F. \ Flipping calculus. \emph{PRIMUS} \ \textbf{2013}, \emph{23}, 477--486. [\href{https://dx.doi.org/10.1080/10511970.2012.757571}{CrossRef}]

\bibitem{sahin2015flipping} Sahin, A.; Cavlazoglu, B.; Zeytuncu, Y.E. \ Flipping a college calculus course: A case study. \emph{J. Educ. Technol. Soc.} \textbf{2015}, \emph{18}, 142--152. [\href{https://www.jstor.org/stable/pdf/jeductechsoci.18.3.142.pdf}{JSTOR}]

\bibitem{sonner2015impact} Sonnert, G.; Sadler, P.M.; Sadler, S.M.; Bressoud, D.M. \ The impact of instructor pedagogy on college calculus students' attitude toward mathematics. \emph{Int. J. Math. Educ. Sci. Technol.} \textbf{2015}, \emph{46}, 370--387. [\href{https://dx.doi.org/10.1080/0020739X.2014.979898}{CrossRef}]

\bibitem{karjanto2019english} Karjanto, N.; Simon, L. \ English-medium instruction Calculus in Confucian-Heritage Culture: Flipping the class or overriding the culture? \emph{Stud. Educ. Eval.} \textbf{2019}, \emph{63}, 122--135. [\href{https://dx.doi.org/10.1016/j.stueduc.2019.07.002}{CrossRef}]

\bibitem{ziegelmeier2015flipped} Ziegelmeier, L.B.; Topaz, C.M. \ Flipped calculus: A study of student performance and perceptions. \emph{PRIMUS} \textbf{2015}, \emph{25}, 847--860. [\href{https://dx.doi.org/10.1080/10511970.2015.1031305}{CrossRef}]

\bibitem{talbert2014inverting} Talbert, R. \ Inverting the linear algebra classroom. \emph{PRIMUS} \textbf{2014}, \emph{24}, 361--374. [\href{https://dx.doi.org/10.1080/10511970.2014.883457}{CrossRef}]

\bibitem{murphy2016student} Murphy, J.; Chang, J.M.; Suaray, K. \ Student performance and attitudes in a collaborative and flipped linear algebra course. \emph{Int. J. Math. Educ. Sci. Technol.} \textbf{2016}, \emph{47}, 653--673. [\href{https://dx.doi.org/10.1080/0020739X.2015.1102979}{CrossRef}]

\bibitem{novak2017flip} Novak, J.; Kensington-Miller, B.; Evans, T. \ Flip or flop? Students' perspectives of a flipped lecture in mathematics. \emph{Int. J. Math. Educ. Sci. Technol.} \textbf{2017}, \emph{48}, 647--658. [\href{https://dx.doi.org/10.1080/0020739X.2016.1267810}{CrossRef}]

\bibitem{karjanto2017flipped} Karjanto, N.; Lee, S.G. \ Flipped classroom in Introductory Linear Algebra by utilizing Computer Algebra System {\sl SageMath} and a free electronic book. \emph{arXiv} \textbf{2017}, arXiv:1705.00739. 

\bibitem{wilson2013flipped} Wilson, S.G. \ The flipped class: A method to address the challenges of an undergraduate statistics course. \emph{Teach. Psychol.} \textbf{2013}, \emph{40}, 193--199. [\href{https://dx.doi.org/10.1177/0098628313487461}{CrossRef}]

\bibitem{butt2014student} Butt, A. \ Student views on the use of a flipped classroom approach: Evidence from Australia. \emph{Bus. Educ. Accredit.} \textbf{2014}, \emph{6}, 33--43. [\href{https://openresearch-repository.anu.edu.au/bitstream/1885/50584/2/01_Butt_Student_views_on_the_use_of_a_2014.pdf}{ANU Repository}]

\bibitem{poffenberger1959factors} Poffenberger, T.; Norton, D. \ Factors in the formation of attitudes toward mathematics. \emph{J. Educ. Res.} \textbf{1959}, \emph{52}, 171--176. [\href{https://dx.doi.org/10.1080/00220671.1959.10882562}{CrossRef}]

\bibitem{hodges2013improving} Hodges, C.B.; Kim, C. \ Improving college students' attitudes toward mathematics. \emph{TechTrends} \textbf{2013}, \emph{57}, 59--66. [\href{https://dx.doi.org/10.1007/s11528-013-0679-4}{CrossRef}]

\bibitem{lo2017toward} Lo, C.K.; Hew, K.F.; Chen, G. \ Toward a set of design principles for mathematics flipped classrooms: A synthesis of research in mathematics education. \emph{Educ. Res. Rev.} \textbf{2017}, \emph{22}, 50--73. [\href{https://dx.doi.org/10.1016/j.edurev.2017.08.002}{CrossRef}]

\bibitem{karjanto2017attitude} Karjanto, N. \ Attitude toward mathematics among the students at Nazarbayev University Foundation Year Programme. \emph{Int. J. Math. Educ. Sci. Technol.} \textbf{2017}, \emph{48}, 849--863. [\href{https://dx.doi.org/10.1080/0020739X.2017.1285060}{CrossRef}] [\href{https://arxiv.org/abs/1609.00480}{arXiv}]

\bibitem{turra2019flipped} Turra, H.; Carrasco, V.; González, C.; Soval, V.; Yáñez, S. \ Flipped classroom experiences and their impact on engineering students' attitudes towards university-level mathematics. \emph{High. Educ. Pedagog.} \textbf{2019}, \emph{4}, 136--155. [\href{https://dx.doi.org/10.1080/23752696.2019.1644963}{CrossRef}]

\bibitem{buil2019flipped} Buil-Fabrega, M.; Martínez, Casanovas, M.; Ruiz-Munzón, N.; Filho, W.L. \ Flipped classroom as an active learning methodology in sustainable development curricula. \emph{Sustainability} \textbf{2019}, \emph{11}, 4577. [\href{https://dx.doi.org/10.3390/su11174577}{CrossRef}] 

\bibitem{rodriguez2020understanding} Rodríguez-Chueca, J.; Molina-García, A.; García-Ar, A.C.; Pérez, J.; Rodríguez, E. \ Understanding sustainability and the circular economy through flipped classroom and challenge-based learning: An innovative experience in engineering education in Spain. \emph{Environ. Educ. Res.} \textbf{2020}, \emph{26}, 238--252. [\href{https://dx.doi.org/10.1080/13504622.2019.1705965}{CrossRef}]

\bibitem{howell2021engaging} Howell, R.A. \ Engaging students in education for sustainable development: The benefits of active learning, reflective practices and flipped classroom pedagogies. \emph{J. Clean. Prod.} \textbf{2021}, \emph{325}, 129318. [\href{https://dx.doi.org/10.1016/j.jclepro.2021.129318}{CrossRef}]

\bibitem{krantz2015how} Krantz, S.G. \ \emph{How to Teach Mathematics}, 3rd ed.; American Mathematical Society: Providence, RI, USA, 2015.

\bibitem{freeman2014active} Freeman, S.; Eddy, S.L.; McDonough, M.; Smith, M.K.; Okoroafor, N.; Jordt, H.; Wenderoth, M.P. \ Active learning increases student performance in science, engineering, and mathematics. \emph{Proc. Natl. Acad. Sci. USA} \textbf{2014}, \emph{111}, 8410--8415. [\href{https://dx.doi.org/10.1073/pnas.1319030111}{CrossRef}]
\end{thebibliography}
\end{document}